\documentclass[11pt,notitlepage,a4paper,reqno,ps]{amsart}
\usepackage{amsmath,amssymb, amsthm}
\usepackage{graphicx}
\usepackage{tikz}

\usepackage{amstext}
\usepackage{latexsym}
\usepackage{amsfonts}
\usepackage{colortbl}

\newcommand{\COND}[3]{
{\it\linespread{1.3}%
\selectfont
\begin{itemize}
\item[\textbf{#1}] #2
\begin{itemize}
#3
\end{itemize}
\end{itemize}
}}

\newcommand{\SOne}[1]{\{u^{#1}>L\}}
\newcommand{\One}[1]{\mbox{\textbf{1}}_{\SOne{#1}}}

\newcommand{\medint}{-\kern  -,375cm\int}

\theoremstyle{plain}
\newtheorem{theorem}{Theorem}[section]

\newtheorem{lemma}[theorem]{Lemma}
\newtheorem{proposition}[theorem]{Proposition}
\theoremstyle{definition}

\theoremstyle{plain}
\newtheorem{remark}[theorem]{Remark}
\theoremstyle{plain}

\numberwithin{equation}{section} \makeatletter

\renewcommand{\p@enumi}{\thesection.}
\makeatother  \allowdisplaybreaks
\title[Local boundedness]{Local boundedness for weak solutions to some quasilinear elliptic systems}
\author[S. Leonardi - F. Leonetti - C. Pignotti - E. Rocha - V. Staicu]{Salvatore Leonardi - Francesco Leonetti - Cristina Pignotti - \\
Eugenio Rocha - Vasile Staicu}
\address{Salvatore Leonardi: DMI - Department of Mathematics and Informatics, University of Catania, Viale A. Doria 6, 95125 Catania, Italy}
\email{leonardi@dmi.unict.it}

\address{ Francesco Leonetti - Cristina Pignotti:  DISIM - Department of Information Engineering, Computer Science and Mathematics, University of L'Aquila, Via Vetoio snc - Coppito, 67100 L'Aquila, Italy}
 \email{leonetti@univaq.it, pignotti@univaq.it}
\address{ Eugenio Rocha - Vasile Staicu:  CIDMA - Center for Research and Development in Mathematics and Applications,
Department of Mathematics, University of Aveiro, 3810-193 Aveiro, Portugal}
 \email{eugenio@ua.pt, vasile@ua.pt}


\begin{document}

%\maketitle
%{\footnotesize
 %% please put the address of the first author
% \centerline{Dipartimento di Matematica ``U. Dini'', Universit\`a di
% Firenze}
% \centerline{Viale Morgagni 67/A,
%50134 - Firenze, Italy}
%} %% Do not forget to end the {\footnotesize by the sign }

\begin{abstract}
We consider quasilinear elliptic systems in divergence form. In general, we cannot expect
that weak solutions are locally bounded because of De Giorgi's counterexample.
Here we assume a condition on the support of off-diagonal coefficients 
that "keeps away" the counterexample and allows us to prove
local boundedness of weak solutions.
\end{abstract}

\maketitle

%\thanks{\textit{Acknowledgements.}}
\section{Introduction}

\noindent
We consider quasilinear elliptic systems in divergence form
\begin{equation}
\label{intro_divergence_form_system}
- div (a(x,u(x)) Du(x)) = 0,
\quad
x \in \Omega,
\end{equation}
\noindent
where $u: \Omega \subset \mathbb{R}^n \to \mathbb{R}^N$ and $a:\Omega \times \mathbb{R}^N \to \mathbb{R}^{N^2n^2}$ is matrix valued with components $a^{\alpha,\beta}_{i,j}(x,y)$ where $i,j \in \{1,...,n\}$ and 
$\alpha,\beta \in \{1,...,N\}$. Note that (\ref{intro_divergence_form_system}) is a system of $N$ equations
\begin{equation}
\left\{
\begin{array}{l}
-\sum\limits_{i=1}^{n} 
D_i
\left(
\sum\limits_{j=1}^{n}
a_{i,j}^{1,1}(x,u) D_j u^1
+
\sum\limits_{j=1}^{n}
a_{i,j}^{1,2}(x,u) D_j u^2
+
...+
\sum\limits_{j=1}^{n}
a_{i,j}^{1,N}(x,u) D_ju^N
\right)
= 0
\\
-\sum\limits_{i=1}^{n} 
D_i
\left(
\sum\limits_{j=1}^{n}
a_{i,j}^{2,1}(x,u) D_j u^1
+
\sum\limits_{j=1}^{n}
a_{i,j}^{2,2}(x,u) D_j u^2
+
...+
\sum\limits_{j=1}^{n}
a_{i,j}^{2,N}(x,u) D_ju^N
\right)
= 0
\\
...................................................................................................
\\
-\sum\limits_{i=1}^{n} 
D_i
\left(
\sum\limits_{j=1}^{n}
a_{i,j}^{N,1}(x,u) D_j u^1
+
\sum\limits_{j=1}^{n}
a_{i,j}^{N,2}(x,u) D_j u^2
+
...+
\sum\limits_{j=1}^{n}
a_{i,j}^{N,N}(x,u) D_ju^N
\right)
= 0
\end{array}
\right.
\label{intro_all_the_equations}
\end{equation}
\noindent
Let us assume that coefficients $a^{\alpha,\beta}_{i,j}(x,y)$ are measurable with respect to $x$, continuous with respect to $y$, bounded and elliptic. When $N=1$, we have only one equation and the celebrated De Giorgi-Nash-Moser theorem forces weak solutions $u \in W^{1,2}(\Omega)$ to be locally bounded and even H\"older continuous, see section 2.1 in \cite{Mingione-regularity}. The result is no longer true, in general, for systems: De Giorgi's counterexample shows that $u(x) = x/|x|^\gamma$ is a weak solution to a particular system (\ref{intro_all_the_equations}) where $\Omega$ is the ball centered at the origin with radius $1$ and $\gamma>1$ is a suitable exponent; it turns out that $u$ cannot be bounded inside $\Omega$ near the origin; see \cite{deg}, section 3 in \cite{Mingione-regularity} and the recent paper \cite{MonSav}; see also \cite{Podio} and \cite{leonardi}.
Now the effort is finding additional restrictions on the coefficients $a^{\alpha,\beta}_{i,j}$ that keep away De Giorgi's counterexample and allow for local boundedness of weak solutions $u$. The easiest case happens when off-diagonal coefficients vanish, that is
\begin{equation}
\label{intro_diagonal_coefficients}
a^{\alpha,\beta}_{i,j}=0
\quad \text{ for }
\quad
\beta \neq \alpha.
\end{equation}
\noindent
In such a case, the $\alpha$ row of the system is 
\begin{equation}
\label{intro_alpha_equation}
-\sum\limits_{i=1}^{n} 
D_i
\left(
\sum\limits_{j=1}^{n}
a_{i,j}^{\alpha,\alpha}(x,u(x)) D_j u^\alpha (x)
\right)
= 0,
\end{equation}
\noindent
so, if $b_{i,j}(x) = a_{i,j}^{\alpha,\alpha}(x,u(x))$ and $v(x) = u^\alpha(x)$, then we are in the linear scalar case
\begin{equation}
\label{intro_scalar_equation}
-\sum\limits_{i=1}^{n} 
D_i
\left(
\sum\limits_{j=1}^{n}
b_{i,j}(x) D_j v (x)
\right)
= 0
\end{equation}
\noindent
and $v$ turns out to be locally bounded and H\"older continuous. A further step has been made in \cite{yan}: the system 
(\ref{intro_all_the_equations}) is assumed to be tridiagonal, that is
\begin{equation}
\label{intro_tridiagonal_coefficients}
a^{\alpha,\beta}_{i,j}=0
\quad 
\text{ for }
\quad
\beta > \alpha.
\end{equation}
In such a case the system (\ref{intro_all_the_equations}) becomes
\begin{equation}
\left\{
\begin{array}{l}
-\sum\limits_{i=1}^{n} 
D_i
\left(
\sum\limits_{j=1}^{n}
a_{i,j}^{1,1}(x,u) D_j u^1
\phantom
{
+
\sum\limits_{j=1}^{n}
a_{i,j}^{1,2}(x,u) D_j u^2
+
...+
\sum\limits_{j=1}^{n}
a_{i,j}^{1,N}(x,u) D_ju^N
}
\right)
= 0
\\
-\sum\limits_{i=1}^{n} 
D_i
\left(
\sum\limits_{j=1}^{n}
a_{i,j}^{2,1}(x,u) D_j u^1
+
\sum\limits_{j=1}^{n}
a_{i,j}^{2,2}(x,u) D_j u^2
\phantom
{
+
...+
\sum\limits_{j=1}^{n}
a_{i,j}^{2,N}(x,u) D_ju^N
}
\right)
= 0
\\
...................................................................................................
\\
-\sum\limits_{i=1}^{n} 
D_i
\left(
\sum\limits_{j=1}^{n}
a_{i,j}^{N,1}(x,u) D_j u^1
+
\sum\limits_{j=1}^{n}
a_{i,j}^{N,2}(x,u) D_j u^2
+
...+
\sum\limits_{j=1}^{n}
a_{i,j}^{N,N}(x,u) D_ju^N
\right)
= 0
\end{array}
\right.
\label{intro_tridiagonal_equations}
\end{equation}
\noindent
Then we can apply to the first equation the regularity for scalar case: 
H\"older continuity for $u^1$ and suitable decay on balls for $Du^1$.
Now the second row can be written as follow
\begin{equation}
-\sum\limits_{i=1}^{n} 
D_i
\left(
\sum\limits_{j=1}^{n}
a_{i,j}^{2,2}(x,u) D_j u^2
\right)
= 
\sum\limits_{i=1}^{n} 
D_i
\left(
\sum\limits_{j=1}^{n}
a_{i,j}^{2,1}(x,u) D_j u^1
\right);
\label{intro_second_equations}
\end{equation}
\noindent
the good behaviour of $Du^1$ on the right hand side can be transferred to the left hand side so that 
$u^2$ inherits H\"older continuity and $Du^2$ gets a suitable decay on balls.
The procedure can be iterated until we arrive at $u^N$.
Another step has been made in \cite{ME} where the local boundedness is obtained 
under the following structure assumption:
there exist numbers $\lambda>0$, $L \geq 0$ and two nonnegative functions $d(x), g(x)$, such that
\begin{equation}
\label{intro_es}
%\begin{array}{lll}
% && 
\sum\limits_{\alpha=1}^{N} \sum\limits_{\gamma=1}^{N}
	\dfrac{y^{\alpha} y^{\gamma}}{|y|^2}   \left( \sum\limits_{i=1}^{n} p^{\gamma}_i \sum\limits_{\beta=1}^{N} \sum\limits_{j=1}^{n} a_{i,j}^{\alpha, \beta}(x,y)p^{\beta}_j \right)  
	%\\
%& 
\geq 
%& 
- \left\{ \delta |p|^2  + \left( \dfrac{1}{\delta}\right)^{\lambda} [d(x) |y|^2 + g(x)] \right\}
%\end{array}
\end{equation}
is fulfilled for all $\delta \in (0,1)$ and all $(x, y, p)$, with $|y|>L$.
In the present work we assume a condition on the support of off-diagonal coefficients:
there exists $L_0\in(0,+\infty)$ such that $\forall \ L\ge L_0,$ when $\alpha\neq\beta,$
\begin{equation}
\label{intro_our_assumtion}
\begin{array}{l}
(a_{i,j}^{\alpha,\beta}\left(  x,y\right)  \neq 0 \ \mbox{and}\  \ y^\alpha > L) \Rightarrow y^\beta>L,
\\
(a_{i,j}^{\alpha,\beta}\left(  x,y\right)  \neq 0 \ \mbox{and}\  y^\alpha < -L) \Rightarrow y^\beta<-L
\end{array}
\end{equation}
(see Figure~\ref{figura_diagonale}). Under such a restriction we are able to prove local boundedness of weak solutions.
All the necessary assumptions and the result will be listed in section 2 while proofs will be performed in section 3.
Let us mention that off-diagonal coefficients with a particular support have been successfully used when proving maximum principles in 
\cite{Leo_Leo_Pign_Roc_Vas} and when obtaining existence for measure data problems in \cite{Leo_Roc_Vas_1}, \cite{Leo_Roc_Vas_2}.
It is worth mentioning that, when the ratio between the largest and the smallest eigenvalues of $a^{\alpha,\beta}_{i,j}$ is close to 1, then regularity of $u$ is studied at page 183 of \cite{giaquinta}; see also \cite{necas}, \cite{kottas}, \cite{koshelev}, \cite{leonardi-kottas-stara}. Let us also say that proving boundedness for weak solutions could be an important tool for getting fractional differentiability, see
the estimate after (4.15) in \cite{esposito-leonetti-mingione}; sometimes, a gain in fractional differentiability can be iterated as in Theorem 3.III of \cite{campanato-cannarsa} and in Theorem 3.3 of \cite{floridia-ragusa}.

\section{Assumptions and Result}

\noindent
Assume
$\Omega$ is an open bounded subset of $\mathbb{R}^{n}$, with $n \geq 2$.
 %For  the sake of brevity, $[k]$ denotes the set $\{1,..,k\}$ when $k \geq 1$ is an integer.
\noindent
Consider the system of $N \geq 2$ equations
\begin{equation}\label{MProb}
-\sum\limits_{i=1}^{n}\frac{\partial}{\partial x_i}\left(\sum\limits_{\beta=1}^{N}\sum\limits_{j=1}^{n}a_{i,j}^{\alpha,\beta}\left(  x,u\right) \frac{\partial}{\partial x_j}u^{\beta}\right)  = 0  \mbox{ in }\Omega, \mbox{ for } \alpha = 1,...,N.
\end{equation}
Note that $u^\beta$ is the $\beta$ component of $u = (u^1,u^2,...,u^N)$. We list our structural conditions.

\begin{figure}[ht!]
 \begin{center}
\begin{tikzpicture}[xscale=0.6,yscale=0.6]
%\draw [dotted] (2,0) -- (2,1);
%\draw [dotted] (2.5,0) -- (2.5,2);
%\draw [dotted] (0,2) -- (1,2);
%\draw [dotted] (0,2.5) -- (2,2.5);
%\draw [dotted] (0,2.75) -- (2.5,2.75);
%\draw [dotted] (0,2.875) -- (2.75,2.875);
%\node at (4.2,4.2) {.};
%\node at (4.3,4.3) {.};
%\node at (4.4,4.4) {.};
%\node at (-4.2,-4.2) {.};
%\node at (-4.3,-4.3) {.};
%\node at (-4.4,-4.4) {.};
\fill[gray!50!white]  (-12,0) rectangle (12,2);
%\fill[gray!50!white]  (2,2) rectangle (12,3);
%\fill[gray!50!white]  (3,3) rectangle (12,3.5);
%\fill[gray!50!white]  (3.5,3.5) rectangle (12,3.75);
%\fill[gray!50!white]  (3.75,3.75) rectangle (12,3.875);
%\fill[gray!50!white] (3.875,3.875) rectangle (12,3.935);
%\fill[gray!50!white] (3.935,3.935) rectangle (12,3.965);
%\fill[gray!50!white] (3.965,3.965) rectangle (12,3.98);
%\fill[gray!50!white] (3.98,3.98) rectangle (12,3.9999);
%\fill[gray!50!white] (4,4) rectangle (12,6);
%\fill[gray!50!white] (6,6) rectangle (12,9);
%\fill[gray!50!white] (9,9) rectangle (12,10.5);
%\fill[gray!50!white] (10.5,10.5) rectangle (12,11.25);
%\fill[gray!50!white] (11.25,11.25) rectangle (12,11.625);
%\fill[gray!50!white] (11.625,11.625) rectangle (12,11.8125);
%\fill[gray!50!white] (11.8125,11.8125) rectangle (12,11.9);
%\fill[gray!50!white] (11.9,11.9) rectangle (12,11.95);
%\fill[gray!50!white] (11.95,11.95) rectangle (12,11.975);
\fill[gray!50!white]  (12,0) rectangle (-12,-2);
%\fill[gray!50!white]  (-2,-2) rectangle (-12,-3);
%\fill[gray!50!white]  (-3,-3) rectangle (-12,-3.5);
%\fill[gray!50!white]  (-3.5,-3.5) rectangle (-12,-3.75);
%\fill[gray!50!white]  (-3.75,-3.75) rectangle (-12,-3.875);
%\fill[gray!50!white] (-3.875,-3.875) rectangle (-12,-3.935);
%\fill[gray!50!white] (-3.935,-3.935) rectangle (-12,-3.965);
%\fill[gray!50!white] (-3.965,-3.965) rectangle (-12,-3.98);
%%\fill[gray!50!white] (-3.98,-3.98) rectangle (-12,-3.9999);
%\fill[gray!50!white] (-4,-4) rectangle (-12,-6);
%\fill[gray!50!white] (-6,-6) rectangle (-12,-9);
%\fill[gray!50!white] (-9,-9) rectangle (-12,-10.5);
%\fill[gray!50!white] (-10.5,-10.5) rectangle (-12,-11.25);
%\fill[gray!50!white] (-11.25,-11.25) rectangle (-12,-11.625);
%\fill[gray!50!white] (-11.625,-11.625) rectangle (-12,-11.8125);
%\fill[gray!50!white] (-11.8125,-11.8125) rectangle (-12,-11.9);
%\fill[gray!50!white] (-11.9,-11.9) rectangle (-12,-11.95);
%\fill[gray!50!white] (-11.95,-11.95) rectangle (-12,-11.975);
\filldraw[gray!50!white] (0,0) -- (12,0) -- (12,12);
\filldraw[gray!50!white] (0,0) -- (-12,0) -- (-12,-12);
%\fill[gray!50!white]  (0,0) rectangle (1,1);
%\fill[gray!50!white]  (0,0) rectangle (2,2);
%\fill[gray!50!white]  (0,0) rectangle (-1,1);
%\fill[gray!50!white] (0,0) rectangle (1,-1);
%\fill[gray!50!white] (0,0) rectangle (-1,-1);
%\fill[gray!50!white]  (1,1) rectangle (2,2);
%\fill[gray!50!white](-1,-1) rectangle (-2,-2);
%\fill[gray!50!white] (-2,-2) rectangle (-3,-3);
%\fill[gray!50!white] (-3,-3) rectangle (-4,-4);
%\draw (2,0.1) -- (2,-0.1);
%\draw (3,0.1) -- (3,-0.1);
%\draw (3.5,0.1) -- (3.5,-0.1);
%\draw (3.75,0.1) -- (3.75,-0.1);
%\draw (3.875,0.1) -- (3.875,-0.1);
%\draw (4,0.1) -- (4,-0.1);
%\draw (-2,0.1) -- (-2,-0.1);
%\draw (-3,0.1) -- (-3,-0.1);
%\draw (-3.5,0.1) -- (-3.5,-0.1);
%\draw (-4,0.1) -- (-4,-0.1);
%\draw (0.1,2) -- (-0.1,2);
%\draw (0.1,3) -- (-0.1,3);
%\draw (0.1,3.5) -- (-0.1,3.5);
%\draw (0.1,3.75) -- (-0.1,3.75);
%\draw (0.1,3.875) -- (-0.1,3.875);
%\draw (0.1,4) -- (-0.1,4);
%\draw (0.1,-2) -- (-0.1,-2);
%\draw (0.1,-3) -- (-0.1,-3);
%\draw (0.1,-3.5) -- (-0.1,-3.5);
%\draw (0.1,-4) -- (-0.1,-4);
%\node [below,scale=0.7] at (2,0) {$k_0^{(1)}$};
%\node [below,scale=0.7] at (3,0) {$k_1^{(1)}$};
%\node [below,scale=0.7] at (3.5,0) {$k_2^{(1)}$};
%\node [below,scale=0.7] at (4,0) {$d^{(1)}$};
%\node [left,scale=0.7] at (0,2) {$k_0^{(1)}$};
%\node [left,scale=0.7] at (0,3) {$k_1^{(1)}$};
%\node [left,scale=0.7] at (0,3.5) {$k_2^{(1)}$};
%\node [left,scale=0.7] at (0,4) {$d^{(1)}$};
\node [left,scale=0.9] at (0,11) {$y^\alpha$};
\node [below,scale=0.9] at (11,0) {$y^\beta$};
\draw [<->] (0,11) -- (0,0) -- (11,0);
\draw (0,-11) -- (0,0) -- (-11,0);
\node [scale=1.2] at (-6,6) {$a^{\alpha,\beta}_{i,j} = 0$};
\node [scale=1.2] at (6,-6) {$a^{\alpha,\beta}_{i,j} = 0$};
\end{tikzpicture}
   \caption{Assumption ($\mathcal{A}_3)$: off-diagonal entries $a_{i,j}^{\alpha, \beta}$ vanish on the white part of the picture; they might be non zero only on the grey part.}
\label{figura_diagonale}
  \end{center}
\end{figure}
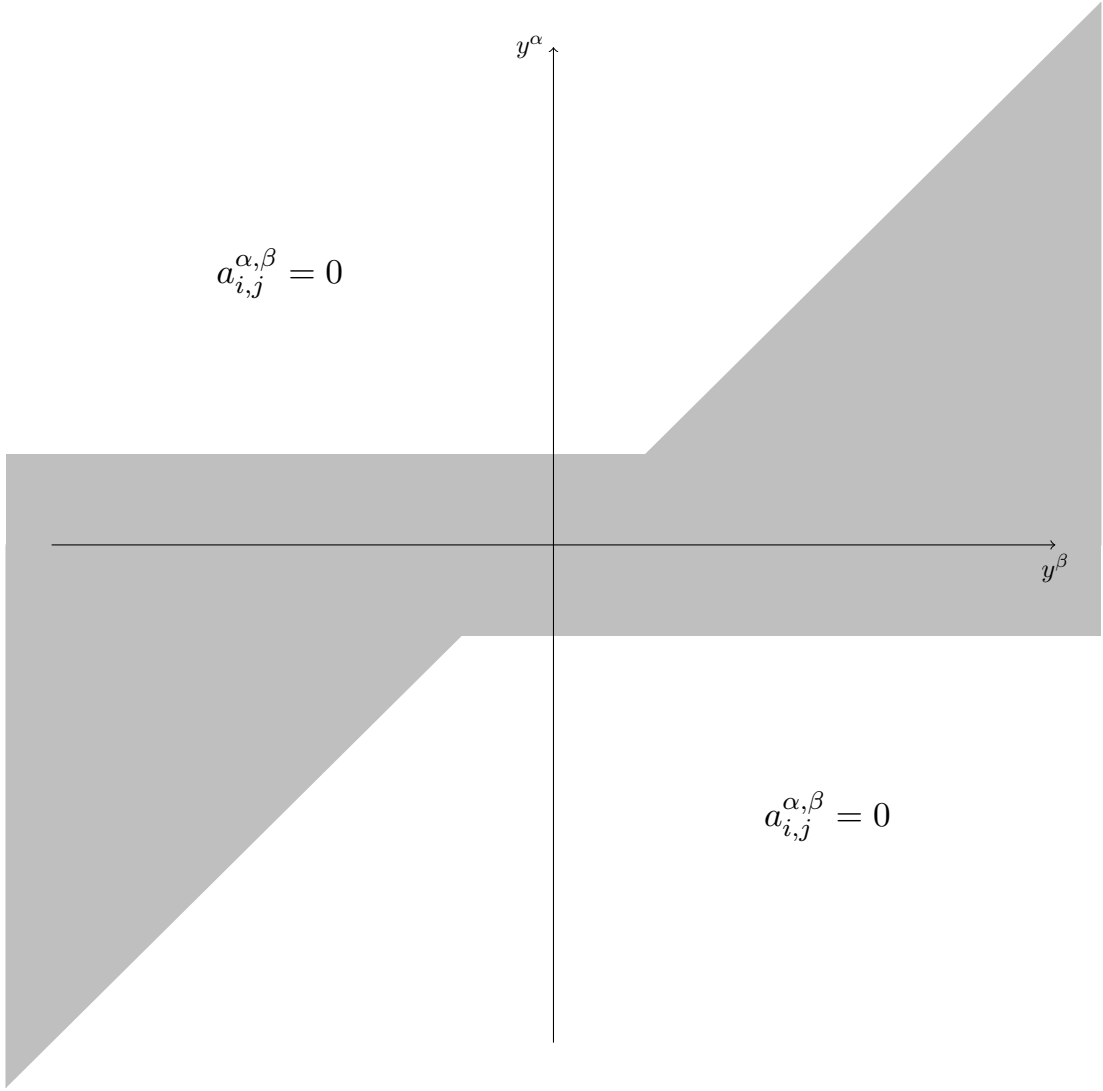

\COND{$(\mathcal{A})$}{For  all $i,j\in\{1,...,n\}$ and all $\alpha,\beta\in\{1,...,N\}$, we require that $a_{i,j}^{\alpha,\beta}:\Omega\times\mathbb{R}^{N}\rightarrow\mathbb{R}$ satisfies the following conditions:}{
\item[$(\mathcal{A}_0)$] $x\mapsto a_{i,j}^{\alpha,\beta}(x,y)$ is measurable and $y\mapsto
a_{i,j}^{\alpha,\beta}(x,y)$ is continuous;

\item[$(\mathcal{A}_1)$] \textit{(boundedness of all the coefficients)} for some positive constant
$c>0$, we have
$$
|a_{i,j}^{\alpha,\beta}\left(  x,y\right)  |\leq c
$$
for almost all $x\in\Omega$ and for all $y\in\mathbb{R}^{N}$;

\item[$(\mathcal{A}_2)$] \textit{(ellipticity of all the coefficients)} for some positive constant
$\nu>0$, we have
$$
\sum\limits_{\alpha,\beta=1}^{N}
\sum\limits_{i,j=1}^{n}
a_{i,j}^{\alpha,\beta}\left(
x,y\right)  \xi^\alpha_{i}\xi^\beta_{j}\geq \nu|\xi|^{2}
%\quad\mbox{ and }\quad b^{\gamma,\gamma}(x,y)\not\equiv 0,
$$
for almost all $x\in\Omega$, for all $y\in\mathbb{R}^{N}$ and for all $\xi\in\mathbb{R}^{N \times n}$;

\item[$(\mathcal{A}_3)$] \textit{(support of off-diagonal coefficients)}
%to be written!!!!!!!!
%there exists $r\in(0,+\infty)$ such that when $\alpha\neq\beta$,
%\begin{eqnarray}
%a_{i,j}^{\alpha,\beta}\left(  x,y\right)  \neq 0 %\Rightarrow
%\mbox{ implies either } y \in \left\{  |y^{\alpha}|  < r,\: |y^{\beta}| < r\right\}
%\mbox{ or }
%\nonumber
%\\
%y\in
%{\displaystyle\bigcup\limits_{h \in \mathbb{Z}}}
%\left\{  hr < y^{\alpha}  < (h+1)r,\:hr < y^{\beta
%} < (h+1)r\right\},
%\nonumber
%\end{eqnarray}
there exists $L_0\in(0,+\infty)$ such that $\forall \ L\ge L_0,$ when $\alpha\neq\beta,$
$$
\begin{array}{l}
(a_{i,j}^{\alpha,\beta}\left(  x,y\right)  \neq 0 \ \mbox{and}\  \ y^\alpha > L) \Rightarrow y^\beta>L\,,\quad\quad (\mathcal{A}_3^\prime)\\
(a_{i,j}^{\alpha,\beta}\left(  x,y\right)  \neq 0 \ \mbox{and}\  y^\alpha < -L) \Rightarrow y^\beta<-L\,.\quad (\mathcal{A}_3^{\prime\prime})
\end{array}
$$
(see Figure~\ref{figura_diagonale}).
}

\noindent
We say that a function $u:\Omega\rightarrow\mathbb{R}^{N}$ is a weak solution of the system (\ref{MProb}), if  $u\in W^{1,2}\left(  \Omega,\mathbb{R}^{N}\right)$ and
\begin{equation}
\label{Sol}
\int_{\Omega}
\sum\limits_{\alpha,\beta=1}^{N}
\sum\limits_{i,j=1}^{n}
a_{i,j}^{\alpha,\beta}(x,u(x))D_{j}u^{\beta}(x)  D_{i}\varphi^{\alpha}(x)dx = 0,
\end{equation}
for all $\varphi\in W_{0}^{1,2}\left(\Omega,\mathbb{R}^{N}\right)$.

\begin{theorem}
\label{Main_theorem} 
Let $u\in W^{1,2}\left(  \Omega,\mathbb{R}^{N}\right)$ be a weak solution of system~(\ref{MProb}) under the set ($\mathcal{A}$) of assumptions. Then $u\in L^{\infty}_{loc}\left(  \Omega,\mathbb{R}^{N}\right)$.
\end{theorem}

\section{Proof of the result}

\noindent
The proof of the Theorem \ref{Main_theorem} will be performed in several steps.

\noindent
\subsection*{STEP 1. Caccioppoli inequality}

\noindent

\begin{theorem}
\label{Caccioppoli}
(Caccioppoli inequality on superlevel sets) 
Let $u\in W^{1,2}\left(  \Omega,\mathbb{R}^{N}\right)$ be a weak solution of system~(\ref{MProb}) under the assumptions ($\mathcal{A}_0$), ($\mathcal{A}_1$), ($\mathcal{A}_2$), ($\mathcal{A}_3^\prime$).
For $0<s<t$, let $B(x_0,s)$ and $B(x_0,t)$ be concentric open balls centered at $x_0$ with radii $s$ and $t$ respectively.
Assume that $B(x_0,t) \subset \Omega$ and $L \geq L_0$. Then
\begin{equation}
\sum\limits_{\alpha=1}^{N}
\int\limits_{\{ u^{\alpha} > L \} \cap  B(x_0,s)}
|D \,u^\alpha|^2  \,dx
\leq
\frac{16 c^2 n^4 N^4}{\nu^2}
\sum\limits_{\alpha=1}^{N}
\int\limits_{\{ u^{\alpha} > L \} \cap  B(x_0,t)}
\left(
\frac{u^{\alpha} - L}{t-s}
\right)^2
\,dx,
\label{Caccioppoli_estimate}
\end{equation}
where $c$ is the constant involved in assumption ($\mathcal{A}_1$), 
$\nu$ is given in ($\mathcal{A}_2$) 
and $L_0$ appears in ($\mathcal{A}_3$).
\end{theorem}
\noindent
{\bf Proof of Theorem \ref{Caccioppoli}}
%\begin{proof}
\noindent
 Let $u\in W^{1,2}\left(  \Omega,\mathbb{R}^{N}\right)$ be a weak solution of system~(\ref{MProb}).
Let $\eta:\mathbb{R}^n \to \mathbb{R}$ be the standard cut-off function such that 
$0 \leq \eta \leq 1$,
$\eta \in C^1_0(B(x_0,t))$, with $B(x_0,t) \subset \Omega$ and $\eta = 1$ in $B(x_0,s)$. Moreover, $|D \eta| \leq 2/(t-s)$ in $\mathbb{R}^n$.
For every level $L \geq L_0$,
%be such that
%$$\sup_{\partial\Omega} u^\alpha \leq L\quad\mbox{ for all } \alpha\in[N].$$
consider the test function $\varphi:\mathbb{R}^n \to \mathbb{R}^N$ with $\varphi = (\varphi^1,...,\varphi^N)$,
where
$$ \varphi^\alpha(x):= \eta^2(x) \max\{0, u^\alpha(x)-L\},\quad\mbox{ for all } \alpha\in\{1,...,N\}. $$
Then
$$  D_i\,\varphi^\alpha=\eta^2 \One{\alpha} D_i\,u^\alpha + 2 \eta (D_i \eta) \One{\alpha} (u^\alpha - L)  \quad\mbox{ for all } i\in\{1,...,n\} \mbox{ and } \alpha\in\{1,...,N\}.$$
Using this  test function in the weak formulation (\ref{Sol}) of system (\ref{MProb}), we have
\begin{align*}
&0 = \int_{\Omega}
\sum\limits_{\alpha,\beta=1}^{N}
\sum\limits_{i,j=1}^{n}
a_{i,j}^{\alpha,\beta}D_{j}u^{\beta}  D_{i}\varphi^{\alpha}\,dx
=
\\
&\int_{\Omega}
\sum\limits_{\alpha,\beta=1}^{N}
\sum\limits_{i,j=1}^{n}
a_{i,j}^{\alpha,\beta}
D_{j}u^{\beta} \eta^2 \One{\alpha} D_i\,u^\alpha\,dx
+
\int_{\Omega}
\sum\limits_{\alpha,\beta=1}^{N}
\sum\limits_{i,j=1}^{n}
a_{i,j}^{\alpha,\beta}
D_{j}u^{\beta}
2 \eta (D_i \eta) \One{\alpha} (u^\alpha - L)dx.
\end{align*}
Now, the assumption ($\mathcal{A}_3^\prime$) guarantees that
\begin{equation}
\label{using_staircase_assumption}
a_{i,j}^{\alpha,\beta}(x,u(x)) \One{\alpha}(x) = a_{i,j}^{\alpha,\beta}(x,u(x)) \One{\beta}(x) \One{\alpha}(x)
\end{equation}
when $\beta \neq \alpha$ and $L\ge L_0.$ It is worthwhile to note that (\ref{using_staircase_assumption}) holds true when $\alpha = \beta$ as well;
then
\begin{eqnarray}
\int_{\Omega}
\sum\limits_{\alpha,\beta=1}^{N}
\sum\limits_{i,j=1}^{n}
a_{i,j}^{\alpha,\beta}
\One{\beta}
D_{j}u^{\beta} \eta^2 \One{\alpha} D_i\,u^\alpha\,dx
\qquad
\qquad
\qquad
\qquad
\nonumber
\\
=
-\int_{\Omega}
\sum\limits_{\alpha,\beta=1}^{N}
\sum\limits_{i,j=1}^{n}
a_{i,j}^{\alpha,\beta}
\One{\beta}
D_{j}u^{\beta}
2 \eta (D_i \eta) \One{\alpha} (u^\alpha - L)\, dx.
%\nonumber
\label{aggiunta_salvatore}
\end{eqnarray}
Now we can use the ellipticity assumption ($\mathcal{A}_2$) with $\xi_i^\alpha = \One{\alpha} D_i\,u^\alpha$ and we get
\begin{eqnarray}
\nu
\int_{\Omega}
\eta^2
\sum\limits_{\alpha=1}^{N}
%\sum\limits_{i=1}^{n}
 \One{\alpha} |D \,u^\alpha|^2  \,dx
\leq
\int_{\Omega}
\sum\limits_{\alpha,\beta=1}^{N}
\sum\limits_{i,j=1}^{n}
a_{i,j}^{\alpha,\beta}
\One{\beta}
D_{j}u^{\beta} \eta^2 \One{\alpha} D_i\,u^\alpha\,dx.
\label{estimate_for_left_hand_side}
\end{eqnarray}
Moreover
\begin{eqnarray}
-\int_{\Omega}
\sum\limits_{\alpha,\beta=1}^{N}
\sum\limits_{i,j=1}^{n}
a_{i,j}^{\alpha,\beta}
\One{\beta}
D_{j}u^{\beta}
2 \eta (D_i \eta) \One{\alpha} (u^\alpha - L)\, dx
\leq
\nonumber
\\
\int_{\Omega}
c
\sum\limits_{\beta=1}^{N}
\sum\limits_{j=1}^{n}
\One{\beta}
|D_{j}u^{\beta}|
\sum\limits_{\alpha=1}^{N}
\sum\limits_{i=1}^{n}
2 \eta |D_i \eta| \One{\alpha} (u^\alpha - L)\,dx
\leq
\nonumber
\\
\int_{\Omega}
c
\sum\limits_{\beta=1}^{N}
n
\One{\beta}
|D u^{\beta}|
\sum\limits_{\alpha=1}^{N}
n
2 \eta |D \eta| \One{\alpha} (u^\alpha - L)\, dx
\leq
\nonumber
\\
\int_{\Omega}
c n^2
\epsilon
\eta^2
\left(
\sum\limits_{\beta=1}^{N}
\One{\beta}
|D u^{\beta}|\right)^2
+
%\nonumber
%\\
\int_{\Omega}
\frac{c n^2}{\epsilon}
|D\eta|^2
\left(
\sum\limits_{\alpha=1}^{N}
\One{\alpha}
(u^{\alpha} - L)
\right)^2\,dx
\leq
\nonumber
\\
\int_{\Omega}
c n^2 N^2
\epsilon
\eta^2
\sum\limits_{\beta=1}^{N}
\One{\beta}
|D u^{\beta}|^2
+
%\nonumber
%\\
\int_{\Omega}
\frac{c n^2N^2}{\epsilon}
|D\eta|^2
\sum\limits_{\alpha=1}^{N}
\One{\alpha}
(u^{\alpha} - L)^2\, dx ,
\label{estimate_for_right_hand_side}
\end{eqnarray}
where we used the inequality $2ab \leq \epsilon a^2 + b^2/\epsilon$, provided $\epsilon > 0$.
Merging (\ref{estimate_for_left_hand_side}) and (\ref{estimate_for_right_hand_side}) into (\ref{aggiunta_salvatore})
we get
\begin{eqnarray}
\nu
\int_{\Omega}
\eta^2
\sum\limits_{\alpha=1}^{N}
%\sum\limits_{i=1}^{n}
 \One{\alpha} |D \,u^\alpha|^2  \,dx
\leq
\nonumber
\\
\int_{\Omega}
c n^2 N^2
\epsilon
\eta^2
\sum\limits_{\beta=1}^{N}
\One{\beta}
|D u^{\beta}|^2
+
%\nonumber
%\\
\int_{\Omega}
\frac{c n^2 N^2}{\epsilon}
|D\eta|^2
\sum\limits_{\alpha=1}^{N}
\One{\alpha}
(u^{\alpha} - L)^2\, dx.
%\label{estimate_for_right_hand_side}
\nonumber
\end{eqnarray}
We take $\epsilon=\nu/(2cn^2N^2)$ and we have
\begin{equation}
\frac{\nu}{2}
\int_{\Omega}
\eta^2
\sum\limits_{\alpha=1}^{N}
 \One{\alpha} |D \,u^\alpha|^2  \,dx
\leq
\int_{\Omega}
\frac{2 c^2 n^4 N^4}{\nu}
|D\eta|^2
\sum\limits_{\alpha=1}^{N}
\One{\alpha}
(u^{\alpha} - L)^2\, dx.
%\label{estimate_for_right_hand_side}
\nonumber
\end{equation}
Using the properties of the cut off function $\eta$ we get
\begin{equation}
\sum\limits_{\alpha=1}^{N}
\int\limits_{\{ u^{\alpha} > L \} \cap  B(x_0,s)}
|D \,u^\alpha|^2  \,dx
\leq
\frac{16 c^2 n^4 N^4}{\nu^2}
\sum\limits_{\alpha=1}^{N}
\int\limits_{\{ u^{\alpha} > L \} \cap  B(x_0,t)}
\left(
\frac{u^{\alpha} - L}{t-s}
\right)^2
\,dx.
\label{final_estimate}
\end{equation}

\noindent
 This ends the proof of Theorem \ref{Caccioppoli}. %\ref{max_min}.
\qed
%\end{proof}
\medskip

%\section{Local boundedness}
%\label{s:boundedness}

%\noindent
%\subsection*{STEP 1. Caccioppoli inequality}

\noindent
\subsection*{The ``excess'' on superlevel sets.}
\noindent
 %STEP 2: Decay of the ``excess'' on superlevel sets.

\noindent
In the previous Caccioppoli inequality the following sum %of integrals 
appears on the right hand side:
\begin{equation}
\sum\limits_{\alpha=1}^{N}
\int\limits_{\{ u^{\alpha} > L \} \cap  B(x_0,t)}
\left(
%\frac{
u^{\alpha} - L
%}{t-s}
\right)^2
\,dx.
\label{excess}
\end{equation}
Note that 
%such  a 
the
sum 
(\ref{excess}) 
%of integrals 
%an "excess" 
%integral 
is zero if and only if all the superlevel sets have zero measure, that is 
$|\{ u^{1} > L \}|= 0$, $|\{ u^{2} > L \}|= 0$, ... , $|\{ u^{N} > L \}|= 0$,
where $|A|$ is the $n$-dimensional Lebesgue measure of $A \subset \mathbb{R}^n$.
This happens when $L \geq \max\{ \text{ esssup } u^1, \text{ esssup } u^2, ... , \text{ esssup } u^N \}$.
On the contrary, if $L < \max\{ \text{ esssup } u^1, \text{ esssup } u^2, ... , \text{ esssup } u^N \}$, then the sum 
%of integrals 
(\ref{excess}) is positive.
Moreover,
\begin{equation}
L \mapsto 
\sum\limits_{\alpha=1}^{N}
\int\limits_{\{ u^{\alpha} > L \} \cap  B(x_0,t)}
\left(
%\frac{
u^{\alpha} - L
%}{t-s}
\right)^2
\,dx
\qquad
\text{ decreases. }
\label{excess_decreases}
\end{equation}
Then, such a sum (\ref{excess}) 
%of integrals 
measures how much $L$ is far from
$\max\{ \text{ esssup } u^1, \text{ esssup } u^2, ... , \text{ esssup } u^N \}$.
Let us call (\ref{excess}) the excess of $u$ with respect to the level $L$, the "excess" for short.
%Note that such an "excess" is zero if and only if all the superlevel sets have zero measure: 
%$|\{ u^{1} > L \}|= 0$, $|\{ u^{2} > L \}|= 0$, ... , $|\{ u^{N} > L \}|= 0$.
We aim to show that the "excess" is zero for a suitable level $L$.
We first show that the "excess" at level $L_2$ can be estimated by means of the a power $\sigma$ of the "excess" at level $L_1$, for a suitable pair of levels $L_2 > L_1$. Then we iterate the procedure.

\noindent
\subsection*{STEP 2. Decay of the ``excess'' on superlevel sets.}

\noindent
 In general, we consider a
vector valued %Sobolev 
function $v:\Omega\subset \mathbb{R}^n\to \mathbb{R}^N$, $n\ge 2$, with $v=(v^1,...,v^N)$ 
and $v \in W_{\rm loc}^{1,p}(\Omega;\mathbb{R}^N)$, $p\ge 1$, where $\Omega$ is an open set in $\mathbb{R}^n$.
%
% In this step  we prove that if
% \eqref{Caccioppoli} holds then
% $u^1$ is locally bounded. We  will do this by proving a  result  for
% {\em scalar} Sobolev functions $v:\Omega\subset \mathbb{R}^n\to \mathbb{R}$, $n\ge 2$; see Proposition \ref{iterazione} below. }
%
%\medbreak
%\noindent
%Let us assume  that $\Omega$ is an open set in $\mathbb{R}^n$ and
%$v$ is a vector valued  function $v \in W_{\rm loc}^{1,p}(\Omega;\mathbb{R}^N)$, $p\ge 1$.
We fix $B_{R_0} = B(x_0,R_0) \Subset \Omega$, with  $R_0<1$ small enough so that
%$|B_{R_0}|<1$ and
\begin{equation}|B_{R_0}|< 1 \quad \text{and}\quad
\sum\limits_{\alpha = 1}^{N}\int_{B_{R_0}}|v^{\alpha}|^{p^*}\,dx<  1,\label{e:leonetti}
\end{equation} 
where $A \Subset \Omega$ means that the closure $\overline{A}$ is a compact set contained in $\Omega$; 
moreover, $p^*=\frac{np}{n-p}$, if $p<n$, and $p^*$ is any $q>p$, else.
For every $R\in (0,R_0]$ we define
the decreasing sequences
%$k\geq k_0 >0$,
$$
\rho_{h}:=\frac{R}{2}+\frac{R}{2^{h+1}}=\frac{R}{2}\left(1+\frac{1}{2^h}\right), \qquad
\bar{\rho}_h:=\frac{\rho_{h}+\rho_{h+1}}{2}=\frac{R}{2}\left(1+\frac{3}{4\cdot 2^h}\right).
$$
Fixed   a  positive constant
%$d$
$d \geq 1$,
define the  increasing sequence of positive real numbers
\begin{equation}
\label{sequence_of_levels}
k_h:= d\left( 1-\frac{1}{2^{h+1}}\right), \,\,h \in \{0,1,2,...\}.
\end{equation}
Moreover,    define the sequence $(J_{h})$,
$$
J_{h}:= \sum\limits_{\alpha = 1}^{N} \int _{A^{\alpha}_{k_{h},\rho_{h}}}(v^{\alpha} - k_h)^{p^*} \, dx,
$$
where $A^{\alpha}_{k,\rho} = \{ v^{\alpha} > k \} \cap B_\rho$ and $B_\rho = B(x_0,\rho)$.
The following result holds.

\begin{proposition}
(Decay of the excess from step h to step h+1) 
\label{iterazione}
Let  $v \in W_{\rm loc}^{1,p}(\Omega;\mathbb{R}^N)$, $p\ge 1$.
Fix $B(x_0,R_0) \Subset \Omega$, with  $R_0<1$ small enough such that  \eqref{e:leonetti} holds.
If there exists
$0\le \vartheta\le 1$ and $c_0>0$ such that, for every  $0<s<t\le R_0$ and for every  $h \in \{0,1,2,...\}$,
\begin{equation}
\label{e:Caccioppoli}
\sum\limits_{\alpha = 1}^{N} \int_{A^{\alpha}_{k_{h},s}}|D v^{\alpha}|^p\,dx
\le
c_0
\sum\limits_{\alpha = 1}^{N}
\left\{
\int_{A^{\alpha}_{k_{h},t}}\left(\frac{v^{\alpha}-k_{h}}{t-s}\right)^{p^*}\,dx+|A^{\alpha}_{k_{h},t}|^{\vartheta}
\right\},
\end{equation}
then,  for every  $R\in (0,R_0]$ and for every  $h \in \{0,1,2,...\}$,
%By \eqref{6formula} it follows
\[J_{h+1}  \le
 c(\vartheta,R)\left(2^{\frac{p^*p^*}{p}}\right)^h (J_{h})^{\vartheta\frac{p^*}{p}},
\]
with the positive constant $c(\vartheta,R)$ independent of $h$.
\end{proposition}

\begin{remark}
We want to stress that the exponent on the right hand side is $p^*$ larger than the exponent $p$ on the left hand side:
this situation has been studied in the scalar case $N=1$ in \cite{moscariello-nania}, \cite{fusco-sbordone}, \cite{cup-leo-mas}.
\end{remark}

\begin{proof}

\noindent
Notice that $(J_h)$ is a decreasing sequence, since the following chain of inequalities holds:
 \begin{eqnarray}
\label{pre3formula}
J_{h+1}
=
\sum\limits_{\alpha = 1}^{N}
\int_{A^{\alpha}_{k_{h+1},\rho_{h+1}}} (v^{\alpha}-k_{h+1})^{p^*}\,dx
\le
\sum\limits_{\alpha = 1}^{N}
\int_{A^{\alpha}_{k_{h+1},\rho_{h}}} (v^{\alpha}-k_{h+1})^{p^*}\,dx
\le
\nonumber
\\
\sum\limits_{\alpha = 1}^{N}
\int_{A^{\alpha}_{k_{h+1},\rho_{h}}} (v^{\alpha}-k_{h})^{p^*}\,dx
\le
\sum\limits_{\alpha = 1}^{N}
\int_{A^{\alpha}_{k_{h},\rho_{h}}} (v^{\alpha}-k_{h})^{p^*}\,dx
=
J_h
%\qquad \forall\, h.
\end{eqnarray}

\noindent
Let us now define a sequence  $(\zeta_{h})$  of  cut-off functions  in $C_0^{1}(B(x_0, \bar{\rho}_{h}))$, such that
$0\le \zeta_h\le 1$,
 \  $\zeta_h\equiv 1$ in $B_{\rho_{h+1}}$,
$|D\zeta_h| \leq {\frac{2^{h+4}}{R}}$.
%Notice that $R_0<1$ implies  $\frac{2^h}{R}>1$.
%\medskip
If we denote $(v-k_{h+1})_+=\max\{v-k_{h+1},0\}$ we get
\begin{align}
\nonumber
J_{h+1}  &
=
\sum\limits_{\alpha = 1}^{N}
\int _{A^{\alpha}_{k_{h+1},\rho_{h+1}}}(v^{\alpha}-k_{h+1})^{p^*} \zeta^{p^*}_{h}\, dx
\leq
\sum\limits_{\alpha = 1}^{N}
\int _{A^{\alpha}_{k_{h+1},\bar{\rho}_h}}(v^{\alpha}-k_{h+1})^{p^*} \zeta^{p^*}_{h}\, dx
\\ &=
\sum\limits_{\alpha = 1}^{N}
\int_{B_{R}}
 (\zeta_h(v^{\alpha}-k_{h+1})_+)^{p^{*}} \,dx.
\label{0formula}
 \end{align}
Sobolev embedding Theorem and  the properties of $\zeta_h$ yield
\begin{align}&\nonumber
 \int_{B_{R}}
 (\zeta_h(v^{\alpha}-k_{h+1})_+)^{p^{*}} \,dx
\\ &
\nonumber
\leq c \left(\int_{B_{R}}
| D(\zeta_h(v^{\alpha}-k_{h+1})_+)|^{p} \,dx\right)^{\frac{p^*}{p}}
\\ \nonumber &\le c \left\{\left(\int_{B_{R}} |Dv^{\alpha}\zeta_h|^{p}\chi_{\{v^{\alpha}>k_{h+1}\}} \,dx\right)^{\frac{1}{p}}+
\left(\int_{B_{R}} |(v^{\alpha}-k_{h+1})_+D\zeta_h|^p\,dx\right)^{\frac{1}{p}}
\right\}^{p^*}
\\ & \le c\left\{\left(\int_{A^{\alpha}_{k_{h+1},\bar{\rho}_{h}}} |Dv^{\alpha}|^{p} \,dx
\right)^{\frac{1}{p}}
+
\left(\left(\frac{2^{h}}{R}\right)^p\int_{A^{\alpha}_{k_{h+1},{\rho}_{h}}}
(v^{\alpha}-k_{h+1})^p\,dx
\right)^{\frac{1}{p}}
\right\}^{p^*}.
\label{1formula}\end{align}
Note that, when $p=n$, we used $|B_R| \leq 1$ (see (\ref{e:leonetti})) in the Sobolev inequality.
Substituting  $t=\rho_h$ and $s=\bar{\rho}_{h}$ in \eqref{e:Caccioppoli} we deduce
\begin{equation}\label{2formula}
\sum\limits_{\alpha = 1}^{N}
\int_{
A^{\alpha}_{k_{h+1},\bar{\rho}_{h}}}
|Dv^{\alpha}|^{p} \,dx\le
c
\sum\limits_{\alpha = 1}^{N}
 \left\{
 \left(\frac{2^{h}}{R}\right)^{p^*}
 \int_{A^{\alpha}_{k_{h+1},\rho_{h}}}
 |v^{\alpha}-k_{h+1}|^{p^*}\,dx+|A^{\alpha}_{k_{h+1},\rho_{h}}|^{\vartheta}\right\}.
\end{equation}
Collecting \eqref{0formula}, \eqref{1formula}, \eqref{2formula}, we obtain
\begin{eqnarray}
\label{3formula}
J_{h+1}  \le
 c\left\{
\sum\limits_{\alpha = 1}^{N}
\left(\frac{2^{h}}{R}\right)^{p^*}
 \int_{A^{\alpha}_{k_{h+1},\rho_{h}}}(v^{\alpha}-k_{h+1})^{p^*}\,dx
+
\sum\limits_{\alpha = 1}^{N}
|A^{\alpha}_{k_{h+1},\rho_{h}}|^{\vartheta}
 +
\right.
\nonumber
\\
\left.
\sum\limits_{\alpha = 1}^{N}
\left(\frac{2^{h}}{R}\right)^p
\int_{A^{\alpha}_{k_{h+1},{\rho}_{h}}}
(v^{\alpha}-k_{h+1})^p\,dx
\right\}^{\frac{p^*}{p}}.
\end{eqnarray}
Since $z^p\le z^{p^*}+1$ for every   $z\ge 0$, then
\[\left(\frac{2^{h}}{R}\right)^p
\int_{A^{\alpha}_{k_{h+1},{\rho}_{h}}}
(v^{\alpha}-k_{h+1})^p\,dx\le \left(\frac{2^{h}}{R}\right)^{p^*}
%\Big\{
\int_{A^{\alpha}_{k_{h+1},{\rho}_{h}}}
(v^{\alpha}-k_{h+1})^{p^*}\,dx+|A^{\alpha}_{k_{h+1},{\rho}_{h}}|
%\Big \}
,\]
so obtaining
\begin{equation}J_{h+1}  \le
 c\left\{
\sum\limits_{\alpha = 1}^{N}
 \left( \frac{2^{h}}{R}\right)^{p^*}
 \int_{A^{\alpha}_{k_{h+1},\rho_{h}}}(v^{\alpha}-k_{h+1})^{p^*}\,dx
+
\sum\limits_{\alpha = 1}^{N}
|A^{\alpha}_{k_{h+1},\rho_{h}}|^{\vartheta}
+
\sum\limits_{\alpha = 1}^{N}
|A^{\alpha}_{k_{h+1},{\rho}_{h}}|\right\}^{\frac{p^*}{p}}.
\label{4formula} \end{equation}
Since
\[
\sum\limits_{\alpha = 1}^{N}
|A^{\alpha}_{k_{h+1},\rho_{h}}|(k_{h+1}-k_h)^{p^*}\le
\sum\limits_{\alpha = 1}^{N}\int_{A^{\alpha}_{k_{h+1},\rho_{h}}}(v^{\alpha}-k_h)^{p^*}\,dx\le J_h,\]
then
\[
|A^{\beta}_{k_{h+1},\rho_{h}}|\le
\sum\limits_{\alpha = 1}^{N}
|A^{\alpha}_{k_{h+1},\rho_{h}}|\le \frac{J_h}{(k_{h+1}-k_h)^{p^*}}=\left(\frac{2^{h+2}}{d}\right)^{p^*}J_h.\]
Taking also into account that (see (\ref{pre3formula}))
% \begin{equation}
% \label{5formula}
\[
\sum\limits_{\alpha = 1}^{N}
\int_{A^{\alpha}_{k_{h+1},\rho_{h}}}(v^{\alpha}-k_{h+1})^{p^*}\,dx
\leq
\sum\limits_{\alpha = 1}^{N}
\int_{A^{\alpha}_{k_{h+1},\rho_{h}}}(v^{\alpha}-k_{h})^{p^*}\,dx
\le J_h,
\]
%\end{equation}
%
%
% \eqref{5formula},
inequality \eqref{4formula} gives
\begin{equation}J_{h+1}  \le c\left\{
 \left( \frac{2^{h}}{R}\right)^{p^*}J_h+\left(\frac{2^{h}}{d}\right)^{\vartheta p^*}(J_h)^{\vartheta}+
 \left(\frac{2^{h}}{d}\right)^{p^*} J_h
  \right\}^{\frac{p^*}{p}}.
  \label{6formula}
\end{equation}
We keep in mind that $J_h$ is decreasing and $k_0 = d/2 > 0$, so
\[
J_h \leq J_0 =
\sum\limits_{\alpha = 1}^{N}
\int_{A^{\alpha}_{k_{0},\rho_{0}}}(v^{\alpha}-k_{0})^{p^*}\,dx
\leq
\sum\limits_{\alpha = 1}^{N}
\int_{A^{\alpha}_{k_{0},\rho_{0}}}(v^{\alpha})^{p^*}\,dx
\leq
\sum\limits_{\alpha = 1}^{N}
\int_{B_R}|v^{\alpha}|^{p^*}\,dx
\leq
1,
\]
where we used (\ref{e:leonetti}).
Since  $J_h\le 1$ for every $h$ and recalling that $d\ge 1> R_0\ge R$, we get
$$
\begin{array}{l}
\displaystyle{
\left( \frac{2^{h}}{R}\right)^{p^*}J_h+\left(\frac{2^{h}}{d}\right)^{\vartheta p^*}(J_h)^{\vartheta}+
 \left(\frac{2^{h}}{d}\right)^{p^*} J_h\le
%  \left\{ \left( \frac{2^{h}}{R}\right)^{p^*}+ \left(\frac{2^{h}}{d}\right)^{\alpha p^*}+ \left(\frac{2^{h}}{d}\right)^{p^*}
%  \right\} J^{\alpha}_h
 \left\{ 2 \frac{2^{hp^*}}{R^{p^*}}+\frac{2^{h\vartheta p^*}}{R^{\vartheta p^*}}\right\} (J_h)^{\vartheta}}\\
 \hspace{3 cm}
 \displaystyle{
 \le \left( \frac{2}{R^{p^*}}+\frac{1}{R^{\vartheta p^*}}\right)2^{hp^*} (J_h)^{\vartheta}.}
 \end{array}
 $$
By \eqref{6formula} it follows
\[J_{h+1}  \le
c\left\{ \left( \frac{2}{R^{p^*}}+\frac{1}{R^{\vartheta p^*}}\right)2^{hp^*} (J_h)^{\vartheta}\right\}^{\frac{p^*}{p}}
\le c(\vartheta,R)\left(2^{\frac{p^*p^*}{p}}\right)^h (J_h)^{\vartheta\frac{p^*}{p}}.
\]
\end{proof}

\medbreak

\noindent
\subsection*{STEP 3. Iteration}
%and  proof of Theorem \ref{t:boundedness}
%}
\noindent

\noindent
We need the following  classical result, see e.g. \cite{giusti}.
\begin{lemma}\label{lemma2}
Let $\gamma >0$ and let $J_h \in [0, +\infty)$ be such that
\begin{equation}
\label{ipetesi_giusti}
J_{h+1} \leq A\,\lambda^h J_h^{1+\gamma}\qquad \forall h\in \mathbb{N}\cup\{0\},
\end{equation}
\noindent
with $A>0$ and $\lambda>1$.
If  $J_{0} \leq A^{-\frac{1}{\gamma}}\lambda^{-\frac{1}{\gamma^2}}$,
then\  
%$J_h \le \lambda^{-\frac{h}{\gamma}}J_{0}$\  and \  
$\lim_{h\to \infty} J_h=0$.
\end{lemma}

\noindent
\subsection*{STEP 4. Conclusion.}

\noindent

\noindent
We have got all we need to give the proof of Theorem \ref{Main_theorem}. 
Fix
$B_{R_0}=B(x_0, R_0) \Subset \Omega$, with  $R_0 < 1$ small enough such that  $|B_{R_0}|< 1$ and  $\int_{B_{R_0}}|u|^{p^*}\,dx   \leq  1$.
From \eqref{Caccioppoli_estimate} we have that,
for every $0<s<t\le R_0$  and every $h$, $u$ satisfies

\begin{equation}
\sum\limits_{\alpha=1}^{N}
\int_{A^\alpha_{k_h,s}}|D \,u^\alpha|^2  \,dx
\leq c_0
\sum\limits_{\alpha=1}^{N}\left\{
\int_{A^\alpha_{k_h,t}}
\left(
\frac{u^{\alpha} - k_h}{t-s}
\right)^{2^*}\, dx+\vert A^\alpha_{k_h,t}\vert
\right \},
\label{final_estimate2}
\end{equation}
 where  $c_0>0$ is  independent of $s,t,h$, provided $L_0 \leq d/2$.
\noindent Therefore  $u$ satisfies \eqref{e:Caccioppoli} of Proposition \ref{iterazione} with
%$\alpha:=\min\{1,1-\frac{qp^*}{p(p^*-q)},1-\frac{rp^*}{q(p^*-r)}\}$ and
%constant $c_0$ depending on $n, N$ and the constant $\nu$ in the ellipticity assumption $(\mathcal A_2).$
$p=2$ and 
$\vartheta = 1$. Then Proposition \ref{iterazione}, applied to $u$, gives
   %By \eqref{6formula} it follows
\begin{equation}
\label{decadimento_per_J_h}
J_{h+1}  \le
 c(R)\left(2^{\frac{2^*2^*}{2}}\right)^h J^{\frac{2^*}{2}},
\end{equation}
\noindent
with the positive constant $c(R)$ independent of $h.$
% \rho_{h}:=\frac{R}{2}+\frac{R}{2^{h+1}}=\frac{R}{2}(1+\frac{1}{2^h}), \qquad
% \bar{\rho}_h:=\frac{\rho_{h}+\rho_{h+1}}{2}=\frac{R}{2}(1+\frac{3}{4\cdot 2^h}).
% $$
% and
% $$
% k_h:= d\left( 1-\frac{1}{2^{h+1}}\right), \,\,h \in \mathbb{N}.
% $$ with \tcb{ $d>0$ }to be chosen later.
%  By Proposition \ref{lemma2} applied with $v$ replaced by $u^1$, we obtain
%  that  for every  $R\in (0,R_0]$,
% %By \eqref{6formula} it follows
% \[J_{h+1}  \le
%  c(\alpha,R)\left(2^{\frac{p^*p^*}{p}}\right)^h J^{\alpha\frac{p^*}{p}}_h.
% \]
%\tcb{By  the previous step we have that  \eqref{decadimento_per_J_h} holds. }
Let us note that

 \[J_{0}:=\sum_{\alpha=1}^N\int_{A^\alpha_{k_0,\rho_0}}\left(u^\alpha -k_0\right)^{2^*}\,dx = \sum_{\alpha=1}^N\int_{A^\alpha_{\frac d2,R}}\left(u^\alpha -\frac d2\right)^{2^*}\,dx\le \sum_{\alpha=1}^N\int_{A^\alpha_{\frac d2,R}}\vert u^\alpha\vert^{2^*}\,dx\]
 and so
 $$ J_0\rightarrow 0 \quad \mbox{when} \ d\rightarrow +\infty\,.$$
 Therefore, we can choose $d >0$ large enough, so that
 \[J_{0} <
 c(R)^{-\frac{1}{\frac {2^*}{2}-1}}\left(2^{\frac{2^*2^*}{2}}
 \right)^{-\frac{1}{(\frac{2^*}{2}-1)^{2}}}.\] Thus, by Lemma \ref{lemma2} we deduce that
  $\lim_{h\to \infty}J_{h}=0$; 
	%and then, 
	since

  $$
  \begin{array}{l}
  \displaystyle{
  J_h=\sum_{\alpha=1}^N\int_{A^\alpha_{k_h,\rho_h}}\left(u^\alpha -k_h\right)^{2^*}\,dx\ge
  \sum_{\alpha=1}^N\int\limits_{\{u^\alpha >k_h\}\cap B_{R/2}}\left(u^\alpha -k_h\right)^{2^*}\,dx
  }\\
  \hspace{2 cm}
  \displaystyle{
  \ge
  \sum_{\alpha=1}^N\int\limits_{\{u^\alpha >d\}\cap B_{R/2}}\left(u^\alpha -k_h\right)^{2^*}\,dx
  \ge \sum_{\alpha=1}^N\int\limits_{\{u^\alpha >d\}\cap B_{R/2}}\left(u^\alpha -d\right)^{2^*}\,dx
  },
  \end{array}
  $$
  we deduce that $\vert \{u^\alpha >d\}\cap B_{R/2}\vert =0,$ namely
  $u\le d$ a.e. in $B_{\frac{R}{2}}$.
We have so proved that $u$ is  locally bounded from above.

%\tcg{HO TOLTO IL BOX (c.v.d.)}
%\qed

%
%\begin{remark}
%Notice that if $n=3$ and $p=2$ then $p^*=6$. Therefore
%\begin{equation}\label{equiv2}q< \frac{p^*(p^{*}-p)p}{(p^*)^2+(p^{*}-p)p} \Leftrightarrow q<\frac{12}{11}, \  r<	
%\frac{p^*(p^*-p)q}{(p^*)^2+(p^{*}-p)q}\Leftrightarrow r<\frac{6q}{p+q}.\end{equation}
%If $q\in \left[1,\frac{12}{11}\right)$,  then we last condition implies $r\in \left[\frac{3}{5},\frac{24}{37}\right)$.  \end{remark}

\noindent
Now, let $\tilde{u} = -u$. Then $\tilde{u} \in W^{1,2}(\Omega; \mathbb{R}^{N})$  and, since $u$ satifies  (\ref{Sol}), then $\tilde{u}$ satisfies
\begin{equation}
\label{Sol2}
0 =  \int_{\Omega}\sum_{\alpha,\beta=1}^{N}\sum_{i,j=1}^{n}
\tilde{a}_{i,j}^{\alpha,\beta}(x, \tilde{u}(x))D_{j}\tilde{u}^{\beta}(x)  D_{i}\varphi^{\alpha}(x)dx
% ,\,\,\,\, \forall \varphi\in W_{0}^{1,2}\left(\Omega,\mathbb{R}^{N}\right)
\end{equation}
for every $\varphi\in W_{0}^{1,2}\left(\Omega,\mathbb{R}^{N}\right)$,
where
\begin{equation}
\label{sol3}
\tilde{a}_{i,j}^{\alpha,\beta}(x, y) : = a_{i,j}^{\alpha,\beta}(x,  - y).
\end{equation}
\noindent
We observe that the new coefficients, defined by (\ref{sol3}), readily satisfy  conditions  ($\mathcal{A}_{0}$),   ($\mathcal{A}_{1}$), ($\mathcal{A}_{2}$).
\noindent
Moreover, if  $\alpha\neq\beta$  the coefficients  $ \tilde{a}_{i,j}^{\alpha,\beta}\left(  x, y\right)$ satisfy ($\mathcal {A}_3^\prime$) provided
$ {a}_{i,j}^{\alpha,\beta}\left(  x, y\right)$  satisfy ($\mathcal {A}_3^{\prime\prime}$). Therefore, we can argue as above on $\tilde u$ obtaining the estimate from below for $u$.
This ends the proof of Theorem \ref{Main_theorem}.
\qed

\section{An example}

%\begin{example}
\noindent
Let us take $N=2$ and $n=3$; 
%$\zeta_1>0$, $\zeta_2>0$, $\theta\in(0,1)$, and a continuous function $b:\mathbb{R}^N\rightarrow\mathbb{R}$ such that $\lim_{|y|\rightarrow\infty} b(y)e^{-|y|}=0$ for $y\in\mathbb{R}^N$, 
we define the matrices $a^{\alpha,\beta}\equiv a^{\alpha,\beta}(y)$ ($\alpha,\beta\in\{1,2\}$ and $y = (y^1,y^2)$) as
$$
a^{1,1}:=\left(
\begin{array}[c]{ccc}
2 & 0 & 0\\
0 & 2 & 0\\
0 & 0 & 1
\end{array}\right),\:\:
a^{1,2}:=\left(
\begin{array}[c]{ccc}
b(y) & 0 & 0\\
0 & 0 & 0\\
0 & 0 & 0
\end{array}\right),\:\:
a^{2,1}:=\left(
\begin{array}[c]{ccc}
0 & w(y) & 0\\
0 & 0 & 0\\
0 & 0 & 0
\end{array}\right),\:\:
a^{2,2}:=\left(
\begin{array}[c]{ccc}
27 & 0 & 0\\
0 & 1 & 0\\
0 & 0 & 1
\end{array}\right),
$$
where  $b:\mathbb{R}^2 \rightarrow \mathbb{R}$ is a  continuous function such that 
$0 \leq b(y^1, y^2) \leq 2$, $b(k,k+1) = 2$ for every integer $k \geq 2$, $b(0,0) = 2$ and the support of $b$ is contained in the grey part of figure (\ref{figura_diagonale}) with $\alpha = 1$ and $\beta = 2$; moreover, 
$w:\mathbb{R}^2 \rightarrow \mathbb{R}$ is a  continuous function such that 
$-10 \leq w(y^1, y^2) \leq 0$, $w(k+1,k) = -10$ for every integer $k \geq 2$, $w(0,0) = -10$ and the support of $w$ is contained in the grey part of figure (\ref{figura_diagonale}) with $\alpha = 2$ and $\beta = 1$.
It easy to check that 
 assumptions~($\mathcal{A}_0$)--($\mathcal{A}_3$) are satisfied with $c=27$ and $\nu = 1$. 
For the convenience of the reader, let us do the calculations for the ellipticity ($\mathcal{A}_2$).
\begin{eqnarray}
\sum\limits_{\alpha,\beta=1}^{N}
\sum\limits_{i,j=1}^{n}
a_{i,j}^{\alpha,\beta}\left(
x,y\right)  \xi^\alpha_{i}\xi^\beta_{j}
=
\qquad
\qquad
\qquad
\qquad
\qquad
\qquad
\qquad
\qquad
\qquad
\qquad
\qquad
\nonumber
\\
\sum\limits_{i,j=1}^{n}
a_{i,j}^{1,1}\left(
x,y\right)  \xi^1_{i}\xi^1_{j}
+
\sum\limits_{i,j=1}^{n}
a_{i,j}^{1,2}\left(
x,y\right)  \xi^1_{i}\xi^2_{j}
+
\sum\limits_{i,j=1}^{n}
a_{i,j}^{2,1}\left(
x,y\right)  \xi^2_{i}\xi^1_{j}
+
\sum\limits_{i,j=1}^{n}
a_{i,j}^{2,2}\left(
x,y\right)  \xi^2_{i}\xi^2_{j}
=
\nonumber
\\
a_{1,1}^{1,1}\left(
x,y\right)  \xi^1_{1}\xi^1_{1}
+
a_{2,2}^{1,1}\left(
x,y\right)  \xi^1_{2}\xi^1_{2}
+
a_{3,3}^{1,1}\left(
x,y\right)  \xi^1_{3}\xi^1_{3}
+
a_{1,1}^{1,2}\left(
x,y\right)  \xi^1_{1}\xi^2_{1}
+
a_{1,2}^{2,1}\left(
x,y\right)  \xi^2_{1}\xi^1_{2}
+
\nonumber
\\
a_{1,1}^{2,2}\left(
x,y\right)  \xi^2_{1}\xi^2_{1}
+
a_{2,2}^{2,2}\left(
x,y\right)  \xi^2_{2}\xi^2_{2}
+
a_{3,3}^{2,2}\left(
x,y\right)  \xi^2_{3}\xi^2_{3}
=
\nonumber
\\
2  \xi^1_{1}\xi^1_{1}
+
2  \xi^1_{2}\xi^1_{2}
+
1  \xi^1_{3}\xi^1_{3}
+
b\left(
y\right)  \xi^1_{1}\xi^2_{1}
+
w\left(
y\right)  \xi^2_{1}\xi^1_{2}
+
27  \xi^2_{1}\xi^2_{1}
+
1  \xi^2_{2}\xi^2_{2}
+
1  \xi^2_{3}\xi^2_{3}
\geq
\nonumber
\\
2  \xi^1_{1}\xi^1_{1}
+
2  \xi^1_{2}\xi^1_{2}
+
1  \xi^1_{3}\xi^1_{3}
-2  |\xi^1_{1}| |\xi^2_{1}|
-10  |\xi^2_{1}| |\xi^1_{2}|
+
27  \xi^2_{1}\xi^2_{1}
+
1  \xi^2_{2}\xi^2_{2}
+
1  \xi^2_{3}\xi^2_{3}
\geq
\nonumber
\\
2  \xi^1_{1}\xi^1_{1}
+
2  \xi^1_{2}\xi^1_{2}
+
1  \xi^1_{3}\xi^1_{3}
-  \xi^1_{1}\xi^1_{1} 
- \xi^2_{1} \xi^2_{1}
- \frac{5}{\epsilon}  \xi^2_{1} \xi^2_{1}
- 5 \epsilon \xi^1_{2} \xi^1_{2}
+
27  \xi^2_{1}\xi^2_{1}
+
1  \xi^2_{2}\xi^2_{2}
+
1  \xi^2_{3}\xi^2_{3}
=
\nonumber
\\
1  \xi^1_{1}\xi^1_{1}
+
2  \xi^1_{2}\xi^1_{2}
+
1  \xi^1_{3}\xi^1_{3}
- \xi^2_{1} \xi^2_{1}
-25  \xi^2_{1} \xi^2_{1}
-1 \xi^1_{2} \xi^1_{2}
+
27  \xi^2_{1}\xi^2_{1}
+
1  \xi^2_{2}\xi^2_{2}
+
1  \xi^2_{3}\xi^2_{3}
=
\nonumber
\\
1  \xi^1_{1}\xi^1_{1}
+
1  \xi^1_{2}\xi^1_{2}
+
1  \xi^1_{3}\xi^1_{3}
+
1  \xi^2_{1}\xi^2_{1}
+
1  \xi^2_{2}\xi^2_{2}
+
1  \xi^2_{3}\xi^2_{3}
= |\xi|^2,
%\nonumber
%\\
\end{eqnarray}
where we used the inequality $2AB \leq A^2 + B^2$ and $2AB \leq \frac{A^2}{\epsilon} + \epsilon B^2$ with $\epsilon = \frac{1}{5}$.

\noindent
On the other hand, this example satisfies neither assumption~(\ref{intro_diagonal_coefficients}) nor assumption~(\ref{intro_tridiagonal_coefficients}) since the two matrices $a^{1,2}$ and $a^{2,1}$ are not zero.
Moreover, (\ref{intro_es}) does not hold true. Indeed, for every integer $k \geq 2$, let us take
\begin{equation}
y^1 = k+1, \quad y^2 = k, \quad p^1_1 = p^1_2 =t > 0, \quad p^2_i = 0 = p^\alpha_3
\end{equation}
and let us compute the left hand side of (\ref{intro_es}); we have
\begin{eqnarray}
I = 
\sum\limits_{\alpha=1}^{N} \sum\limits_{\gamma=1}^{N}
	\dfrac{y^{\alpha} y^{\gamma}}{|y|^2}   \left( \sum\limits_{i=1}^{n} p^{\gamma}_i \sum\limits_{\beta=1}^{N} \sum\limits_{j=1}^{n} a_{i,j}^{\alpha, \beta}(x,y)p^{\beta}_j \right) =
	\nonumber
	\\
\sum\limits_{\alpha=1}^{2} \sum\limits_{\gamma=1}^{2}
	\dfrac{y^{\alpha} y^{\gamma}}{|y|^2}   \left( \sum\limits_{i=1}^{3} p^{\gamma}_i \sum\limits_{\beta=1}^{2} 
\left(	
a_{i,1}^{\alpha, \beta}(y)p^{\beta}_1 + a_{i,2}^{\alpha, \beta}(y)p^{\beta}_2
\right)
\right) =	
\nonumber
\\
\sum\limits_{\alpha=1}^{2} \sum\limits_{\gamma=1}^{2}
	\dfrac{y^{\alpha} y^{\gamma}}{|y|^2}   \left( \sum\limits_{i=1}^{3} p^{\gamma}_i  
\left(	
a_{i,1}^{\alpha, 1}(y)p^{1}_1 + a_{i,2}^{\alpha, 1}(y)p^{1}_2
\right)
\right) =	
\nonumber
\\
\sum\limits_{\alpha=1}^{2} \sum\limits_{\gamma=1}^{2}
	\dfrac{y^{\alpha} y^{\gamma}}{|y|^2}   \left( 
	p^{\gamma}_1  
\left(	
a_{1,1}^{\alpha, 1}(y)p^{1}_1 + a_{1,2}^{\alpha, 1}(y)p^{1}_2
\right)
+
p^{\gamma}_2  
\left(	
a_{2,1}^{\alpha, 1}(y)p^{1}_1 + a_{2,2}^{\alpha, 1}(y)p^{1}_2
\right)
\right) =	
\nonumber
\\
\sum\limits_{\alpha=1}^{2} 
	\dfrac{y^{\alpha} y^{1}}{|y|^2}   \left( 
	p^{1}_1  
\left(	
a_{1,1}^{\alpha, 1}(y)p^{1}_1 + a_{1,2}^{\alpha, 1}(y)p^{1}_2
\right)
+
p^{1}_2  
\left(	
a_{2,1}^{\alpha, 1}(y)p^{1}_1 + a_{2,2}^{\alpha, 1}(y)p^{1}_2
\right)
\right) =	
\nonumber
\\ 
	\dfrac{y^{1} y^{1}}{|y|^2}   \left( 
	p^{1}_1  
\left(	
a_{1,1}^{1, 1}(y)p^{1}_1 + a_{1,2}^{1, 1}(y)p^{1}_2
\right)
+
p^{1}_2  
\left(	
a_{2,1}^{1, 1}(y)p^{1}_1 + a_{2,2}^{1, 1}(y)p^{1}_2
\right)
\right) +	
\nonumber
\\ 
	\dfrac{y^{2} y^{1}}{|y|^2}   \left( 
	p^{1}_1  
\left(	
a_{1,1}^{2, 1}(y)p^{1}_1 + a_{1,2}^{2, 1}(y)p^{1}_2
\right)
+
p^{1}_2  
\left(	
a_{2,1}^{2, 1}(y)p^{1}_1 + a_{2,2}^{2, 1}(y)p^{1}_2
\right)
\right) =	
\nonumber
\\
\dfrac{y^{1} y^{1}|t|^2}{|y|^2}   \left(   
\left(	
a_{1,1}^{1, 1}(y) + a_{1,2}^{1, 1}(y)
\right)
+  
\left(	
a_{2,1}^{1, 1}(y) + a_{2,2}^{1, 1}(y)
\right)
\right) +	
\nonumber
\\ 
	\dfrac{y^{2} y^{1}|t|^2}{|y|^2}   \left(   
\left(	
a_{1,1}^{2, 1}(y) + a_{1,2}^{2, 1}(y)
\right)
+  
\left(	
a_{2,1}^{2, 1}(y) + a_{2,2}^{2, 1}(y)
\right)
\right) =	
\nonumber
\\
\dfrac{(k+1) (k+1)|t|^2}{(k+1)^2 + k^2}   4 +	
	\dfrac{k (k+1)|t|^2}{(k+1)^2 + k^2}  (-10) =	
\nonumber
\\
\dfrac{(-6k+4) (k+1)|t|^2}{(k+1)^2 + k^2} =
\nonumber
\\
\dfrac{(-6k^2 - 2k +4)|t|^2}{2k^2 + 2k +1}.   	
\nonumber
\end{eqnarray}
Now we compute the right hand side of (\ref{intro_es}); we have
\begin{eqnarray}
Q = 
- \left\{ \delta |p|^2  + \left( \dfrac{1}{\delta}\right)^{\lambda} [d(x) |y|^2 + g(x)] \right\}
=
- \left\{ \delta 2 |t|^2  + \left( \dfrac{1}{\delta}\right)^{\lambda} [d(x) [(k+1)^2 + k^2] + g(x)] \right\}
=
\nonumber
\\
- \delta 2 |t|^2
\left\{ 1 +  \dfrac{1}{2 |t|^2 \delta^{1 + \lambda}} [d(x) [(k+1)^2 + k^2] + g(x)] \right\}.
\nonumber
\end{eqnarray}
Let us take 
\begin{equation}
|t|^2 = \dfrac{5}{2  \delta^{1 + \lambda}} [(d(x) + 1) [(k+1)^2 + k^2] + g(x)]
\nonumber
\end{equation}
so that, since $\delta \in (0,1)$,
\begin{equation}
\delta 2 |t|^2
\left\{ 1 +  \dfrac{1}{2 |t|^2 \delta^{1 + \lambda}} [d(x) [(k+1)^2 + k^2] + g(x)] \right\} \leq 2 |t|^2 
\left(
1+\dfrac{1}{5}
\right)
= \dfrac{12}{5} |t|^2
\nonumber
\end{equation}
and
\begin{equation}
\dfrac{-12}{5} |t|^2
\leq
- \delta 2 |t|^2
\left\{ 1 +  \dfrac{1}{2 |t|^2 \delta^{1 + \lambda}} [d(x) [(k+1)^2 + k^2] + g(x)] \right\} 
= Q.
\nonumber
\end{equation}
For every $L>0$, we take $k$ so large that $|y|^2 = (k+1)^2 + k^2 > L^2$ and 
\begin{equation}
\dfrac{-6k^2 - 2k +4}{2k^2 + 2k +1} < \dfrac{-12}{5}.
\nonumber
\end{equation}
Then,
\begin{eqnarray}
I = \dfrac{(-6k^2 - 2k +4)|t|^2}{2k^2 + 2k +1} <  \dfrac{-12}{5} |t|^2 
\leq Q
\nonumber 
\end{eqnarray}
and this shows that the example does not satisfy (\ref{intro_es}).

%\end{example}

\ \\
\noindent
{\bf Acknowledgement.} 
Leonardi acknowledges support from Piano Triennale Ricerca (UNICT) 2016-2018--linea di intervento 2: "Metodi variazionali ed equazioni differenziali”.
Leonardi, Leonetti, Pignotti have
been partially supported by the Gruppo Nazionale per l'Analisi Matematica,
la Probabilit\`a e le loro Applicazioni (GNAMPA) of the Istituto Nazionale di Alta Matematica (INdAM). 
Leonetti and Pignotti acknowledge the support from RIA-UNIVAQ. 
Rocha and Staicu acknowledge the partial support by the Portuguese Foundation
for Science and Technology (FCT), through CIDMA - Center for Research and
Development in Mathematics and Applications, within project UID/MAT/04106/2019(CIDMA).

\end{document}